\documentclass [a4paper,twoside,10pt]{article}

\usepackage[utf8]{inputenc}
\usepackage{amsmath,amsfonts,amssymb}
\usepackage{vmargin,graphicx,theorem}

\usepackage[english]{babel}
\usepackage{enumerate}
\usepackage{color}
\usepackage{pst-fill,pst-grad,pst-plot,pst-eucl,pstricks-add,pst-node}
\RequirePackage[colorlinks,linkcolor=blue,citecolor=blue,urlcolor=blue]{hyperref}

\setpapersize[portrait]{A4}
\setmarginsrb{1.5cm}{0.7cm}{1.5cm}{2cm}{1.5cm}{0cm}{0.5cm}{2cm}
\newcommand{\disp}{\displaystyle}
\newcommand{\dR}{\ensuremath{\mathbf{R}}}

\newtheorem{ethm}{Theorem}[section]

\newtheorem{edefi}[ethm]{Definition}

\newtheorem{erem}[ethm]{Remark}

\newcommand{\proofend}{~$\rhd$}
\newcommand{\proofbegin}{~$\lhd$}

\newcommand{\PAR}[1]{\ensuremath{{\left(#1\right)}}} % (1)
\newcommand{\SBRA}[1]{\ensuremath{{\left[#1\right]}}} % [1]
\renewcommand{\phi}{\varphi}
\renewcommand{\geq}{\geqslant}

\def\disp{\displaystyle}

\newcommand{\R}{\dR}

\newcommand{\grad}{{\rm grad}}

\newcommand{\beq}{\begin{equation}}\newcommand{\eeq}{\end{equation}}
\usepackage{todonotes}

\begin{document}

\title{On the variational interpretation of local logarithmic Sobolev inequalities}

\author{Gauthier Clerc\thanks{Institut Camille Jordan, Umr Cnrs 5208, Universit\'e Claude Bernard Lyon 1, 43 boulevard du 11 novembre 1918, F-69622 Villeurbanne cedex. clerc@math.univ-lyon1.fr, gentil@math.univ-lyon1.fr},  Giovanni Conforti\thanks{D\'epartement de Math\'ematiques Appliqu\'ees, \'Ecole Polytechnique, Route de Saclay, 91128, Palaiseau Cedex, France. giovanni.conforti@polytechnique.edu},  Ivan Gentil$^*$}

\date{\today}
\maketitle
\begin{center}
	Résumé
\end{center}

Le calcul d'Otto est un outil puissant  pour quantifier la dissipation d'énergie mais surtout il propose une interprétation géométrique simple de certaines inégalités fonctionnelles comme l'inégalité de Sobolev logarithmique. La version {\it locale} de ces inégalités, dont la démonstration repose sur le $\Gamma$-calcul développé par Bakry, \'Emery et Ledoux, n'avait pas encore une interprétation géométrique à la Otto. Dans cette courte note, nous comblons cette lacune et montrons comment le calcul d'Otto appliqué au problème de Schrödinger permet une interprétation variationnelle des inégalités de type Sobolev logarithmique locale, ce qui pourrait être un espoir pour l'élaboration de nouvelles d'inégalités locales. 

\begin{center}
	Abstract 
\end{center}

The celebrated Otto calculus has established itself as a powerful tool for proving quantitative energy dissipation estimates and provides with an elegant geometric interpretation of certain functional inequalities such as the Logarithmic Sobolev inequality. However, the \emph{local} versions of such inequalities, which can be proven by means of Bakry-\'Emery-Ledoux  $\Gamma$-calculus, has not yet been given an interpretation in terms of this Riemannian formalism. In this short note we close this gap by explaining  how Otto calculus applied to the Schrödinger problem yields a variations interpretation of the local logarithmic Sobolev inequalities, that could possibly unlock novel class of local inequalities.

\section{Introduction}

\subsection*{Curvature dimension conditions}

We are working on a $(N,\mathfrak g)$, a smooth, connected  and compact Riemannian manifold without boundary.  The heat equation, $\partial_t u=\Delta u$, starting from $f$ is denoted $(t,x)\mapsto P_tf(x)$, where 
$$
P_tf(x)=\int f(y)p_t^x(y)dy,\,\,x\in N,
$$
and $p_t^x(y)$ is the density of the heat kernel with respect to the Riemannian measure $dy$. $\Delta=\nabla\cdot\nabla$ is the Laplace-Beltrami operator ($\nabla\cdot$ is the divergence operator whereas $\nabla$ is the usual gradient operator). $(P_t)_{t\geq0}$ is a reversible Markov semigroup in the sense of~\cite{bgl-book}. Next, we introduce the carré du champ operator $\Gamma$ and its iterated form $\Gamma_2$, they are defined by 
$$
\Gamma(f)=|\nabla f|^2, \quad \Gamma_2(f)=\frac{1}{2}\Delta\Gamma(f)-\Gamma(f,\Delta f).
$$
In particular, B\"ochner's formula yields the equivalent expression of $\Gamma_2$ in terms of Ricci's curvature:
$$
\Gamma_2(f)=||\nabla f||^2_{\text{H.S.}}+Ric_{\mathfrak{g}}(\nabla f).
$$

Following the seminal paper~\cite{bakry-emery1985}, we introduce for  $\rho\in\R$, $n\in [1,\infty]$ the curvature-dimension condition $CD(\rho,n)$ by asking that 
$$
\Gamma_2(f)\geq \rho\Gamma(f)+\frac{1}{n}(\Delta f)^2,
$$ 
 holds for any smooth function $f$.
In view of B\"ochner's formula, the $CD(\rho,n)$ condition is equivalent to the fact that $Ric_{\mathfrak{g}}\geq \rho$ and $d\leq n$ where $d$ is the dimension of $N$.

We also need to introduce the dual semigroup $(P_t^*)_{t\geq 0}$ acting on absolutely continuous measures $\mu\in\mathcal P(N)$ as follows
\begin{equation}
\label{eq-1}
P_t^*(\mu)=P_t\PAR{\frac{d\mu}{dx}}\,dx.
\end{equation}

For the sake of clarity, we limit the discussion in the present note to the Laplace-Beltrami operator in a compact Riemannian Manifold instead a general diffusion Markov operator. However, all of the results we are going to present, can be formulated in such general context  following the paradigm of~\cite{bgl-book}.

\subsection*{Bakry-\'Emery-Ledoux calculus and local logarithmic Sobolev inequalities}

Bakry-\'Emery-Ledoux calculus, known as $\Gamma$ calculus, has been used as a powerful tool to analyse the quantitative behaviour of the heat semigroup $(P_t)_{t\geq0}$ in various ways. Among its main achievements are \emph{local} logarithmic Sobolev inequalities, proven under $CD(\rho,\infty)$ and the $CD(0,n)$ conditions\footnote{Extension to the general $CD(\rho,n)$ condition can be found in~\cite{bbg17}. However, since $CD(\rho,\infty)$ and $CD(0,n)$ are the most relevant settings in applications, we limit the present discussion to these two scenarios.}.

\begin{enumerate}[\bf 1.]
\item {\bf Local inequalities under $CD(\rho ,\infty )$, cf.~\cite{bakry97}.}  Under the  $CD(\rho ,\infty )$ condition the following inequalities are valid for all positive functions $g$ and $T>0$:
\begin{itemize}
\item The mild gradient commutation estimate
\beq 
\label {eq-10}
\frac{\Gamma(P_{T}g)}{P_Tg}  \leq e^{-2\rho T} P_T \Big(\frac{\Gamma (g)}{g}\Big).
\eeq
\item The local logarithmic Sobolev inequality
\beq  \label {eq-20}
 P_T(g\log g)- P_Tg \log P_T g \leq \frac{1 - e^{-2\rho T}}{2\rho} \, P_T \Big (\frac{\Gamma(g)}{g} \Big).
 \eeq
\item The reverse local logarithmic Sobolev inequality
\beq \label {eq-30}
 P_T(g\log g) - P_Tg\log P_T g \geq \frac{e^{2\rho T}-1}{2\rho} \, \frac{\Gamma(P_Tg)}{P_Tg}. 
\eeq
\end{itemize}
\item {\bf Local inequality under the $CD(0,n)$ condition, cf.~\cite{bakry-ledoux}.} Under the  $CD(\rho ,\infty )$ condition the following inequalities are valid for all positive functions $g$ and $T>0$:
\begin{itemize}
\item The local logarithmic Sobolev inequality
\beq\label {eq-40}
\frac{n}{2T}P_Tg\exp\Big[\frac{2}{nP_Tg}\big( P_T(g\log g) - P_Tg\log P_Tg-TLP_Tg\big)\Big]\leq P_T\left(\frac{\Gamma(g)}{g}\right)-\Delta P_Tg+\frac{n}{2T}P_Tg,
 \eeq
and in particular
\beq\label {eq-60}
\frac{\Delta P_Tg}{P_Tg}-\frac{1}{P_T g}P_T\left(\frac{\Gamma(g)}{g}\right)\leq \frac{n}{2T}
 \eeq
 \item The reverse local logarithmic Sobolev inequality
\beq\label {eq-50}
 \frac{n}{2T}P_Tg\exp\Big[-\frac{2}{nP_Tg}\big( P_T(g\log g) - P_Tg\log P_Tg-TLP_Tg\big)\Big]\leq-\frac{\Gamma(P_Tg)}{P_Tg}+\Delta P_Tg+\frac{n}{2T}P_Tg,
 \eeq
 which yields as a by product the celebrated Li-Yau inequality \cite{li1986parabolic},\cite{bakry-ledoux}.
\beq\label {eq-70}
\frac{\Gamma(P_Tg)}{(P_Tg)^2}-\frac{\Delta P_Tg}{P_Tg}\leq \frac{n}{2T}.
 \eeq
\end{itemize}
\end{enumerate}

\begin{erem}
\begin{itemize}
\item The main interest of such inequalities is that there are local, in the sense that they give estimates on the heat kernel itself instead of the entropy of the heat kernel. Moreover, in contrast with their integrated versions, these inequalities imply back the corresponding curvature-dimension conditions, i.e. inequalities~\eqref{eq-10}-\eqref{eq-30} imply $CD(\rho,\infty)$ and~\eqref{eq-40},\eqref{eq-50} imply $CD(0,n)$. 

%\item Inequality~\eqref{eq-10} is known as the gradient bound (or mild commutation between $\Gamma$ and $P_t$), inequalities~\eqref{eq-20} and~\eqref{eq-40} are called local Logarithmic Sobolev inequalities whereas inequalities~\eqref{eq-30}  and~\eqref{eq-50} are called reverse local Logarithmic Sobolev inequalities.  The last inequality~\eqref{eq-70} is the so-called Li-Yau inequality, equivalent to the parabolic Harnack inequality.  

\item Inequalities~\eqref{eq-10}-\eqref{eq-70} are also proved in~\cite[Theorem~5.5.2 and Theorem~6.7.3]{bgl-book}.
\item Many others inequalities have been proved using the same kind of method. Unfortunately, it is actually not possible to reach it with Otto Calculus and the Schrödinger problem. Let mention for instance the isoperimetric inequality under the $CD(\rho,\infty)$ condition proved in~\cite{bakry-ledoux96}, the local hypercontractivity inequality from~\cite{bakry-bolley2012} etc.  
\end{itemize}
\end{erem}

\subsection*{Otto calculus}

The term Otto calculus, introduced by C. Villani in his book~\cite{villani2009}, refers to a Riemannian formalism on the space of probability measures  $M=\mathcal{P}(N)$ that allows, among other things, to interpret the Wasserstein distance as a geodesic distance and the heat equation as the gradient flow of the entropy. Developed in the seminal papers~\cite{jko1998,otto2001, otto-villani2000}, it has led to major advances in the quantifying trend to equilibrium of possibly non  linear PDEs, to the discovery of new functional inequalities as well as to a new geometric interpretation of classical ones. We briefly present here the main ideas that lead to the  definition of the Riemannian metric known as Otto metric. Our heuristic presentation is based on~\cite{gentil-leonard2020}, to which we refer for more details (see also~\cite{gentil} for an informal  presentation in French). However, given its importance, there are many other references that provide with an introduction to Otto calculus, see e.g. the monograph~\cite{gigli2012} for rigorous constructions,~\cite{ambrosio2013user} for a user-friendly presentation and the book~\cite[Chap.~1]{ambrosio-gigli2008}.

We begin by introducing the two main functionals of interest for this work: the Wasserstein distance and the relative entropy. The first one is defined for all $\mu,\nu\in\mathcal P_2(N)$ by
$$
W_2(\mu,\nu)=\inf \sqrt{\iint d(x,y)^2d\pi(x,y)},
$$
where the infimum  runs over all $\pi\in\mathcal P (N\times N)$ with marginals $\mu$ and $\nu$ and $P_2(N)$ is the set of probability measures with finite second moment (since in our case $N$ is compact, $\mathcal P_2(N)=\mathcal P(N)$). For any probability measure $\mu\in\mathcal P(N)$ the entropy with respect to the Riemannian measure is given by
\begin{equation}
\label{eq-2}
\mathcal{F}(\mu)=\begin{cases} \int \log \left( \frac{d \mu}{dx}\right)\,d\mu, \quad & \mbox{if $d\mu\ll dx$,} \\ 
+\infty, \quad & \mbox{otherwise.} \end{cases}
\end{equation}
We recall that a path $[0,1]\ni t\mapsto \mu_t\in\mathcal P_2(N)$ is absolutely continuous if there exists a non negative function $l\in L^2([0,1])$ such that for any $0\leq s\leq t\leq 1$, 
$$
W_2(\mu_t,\mu_s)\leq\int_s^t l(r)dr.
$$
In that case, one can define 
$$
|\dot\mu_t|:=\underset{s\rightarrow t}{\limsup}\frac{W_2(\mu_t,\mu_s)}{|t-s|}\in L^2([0,1]).
$$
Moreover, for any absolutely continuous path $(\mu_t)_{t\geq0}$, there exists a unique vector field $(t,x) \mapsto V_t(x)$ such that $\int|V_t|^2d\mu_t<\infty$ 
and $ |\dot\mu_t|^2=\int|V_t|^2d\mu_t$ (a.e. in $[0,1]$). The vector field $V_t$ is in fact  a limit in $L^2(\mu_t)$ of gradient of smooth  functions in  $N$. For every $t \in [0,1]$ we set
\begin{equation}\label{eq-100}
\dot \mu_t:=V_t,
\end{equation}
and we call $\dot \mu_t$ the velocity of the path $(\mu_t)_{t \in [0,1]}$ at time $t$. 
Finally, the vector field $V_t$ is a weak solution of the continuity equation
\begin{equation}\label{eq-101}
\partial_t\mu_t+\nabla\cdot(\mu_t V_t)=0.
\end{equation}
To give an example, in the case of the heat equation $\partial_t \nu_t=\Delta\nu_t$, the velocity of the path $(\nu_t)_{t\geq0}$ is 
\begin{equation}
\label{eq-102}
\dot\nu_t=-\nabla\log \frac{d\nu_t}{dx}. 
\end{equation}
The tangent space at $\mu\in\mathcal P_2(N)$ is identified with
\[T_{\mu} \mathcal P_2(N)=\overline{\{\nabla \varphi,\,\,\phi:N \mapsto\R,  \varphi \in C^{\infty}(N)\}}^{L_2(\mu)}.\]
The Riemannian metric on $T_{\mu} \mathcal P_2(N)$ is then defined via the scalar product $L^2(\mu)$, 
$$
\langle\nabla \varphi,\nabla\psi\rangle_\mu=\int \Gamma(\varphi,\psi)d\mu,\,\,\,\nabla \varphi,\nabla\psi\in T_{\mu} \mathcal P_2(N).
$$ 

Such metric is often referred to the Otto metric and it can be seen that the geodesics associated to the Otto metric are the displacement interpolations of optimal transport. Using this, a straightforward computation implies that the gradient of the entropy $\mathcal{F}$ at $\mu$ is given by
$$
\grad_{\mu} \mathcal{F}= \nabla \log \left( \frac{d \mu}{dx}\right) \in T_{\mu} \mathcal P_2(N).
$$
Accordingly, we can rewrite the Fisher information functional $\mathcal{I}$ as
\begin{equation}
\label{eq-104}
\mathcal{I}(\mu):=|\grad_{\mu} \mathcal F|^2_{\mu}=:\mathbf{\Gamma}(\mathcal{F})(\mu),
\end{equation}
where $\mathbf{\Gamma}(\mathcal{F})$ can be interpreted as the carré du champ operator applied to the functional $\mathcal{F}$.

In light of~\eqref{eq-102}, we can now view the semigroup $(P_t^*)_{t\geq 0}$ as the gradient flow of the function $\mathcal{F}$, that is to say 
$$
\dot\nu_t=-\grad_{\nu_t} \mathcal{F}.
$$
Now, we turn our attention to the second order calculus introducing covariant derivatives and Hessians. A remarkable fact is that the Hessian of the entropy $\mathcal{F}$ can be expressed in terms of the $\Gamma_2$ operator. In fact, we have (see for instance~\cite{otto-villani2000} and~\cite[Sec 3.3]{gentil-leonard2020}), %for every $\mu \in \mathcal P_2(N)$ it is a bi-linear form in $T_{\mu} \mathcal P_2(N)$ given by 
$$
\forall\,\mu \in \mathcal P_2(N),\, \nabla \phi, \nabla \psi \in T_{\mu} \mathcal{P}_2(N), \quad\mathrm{Hess}_{\mu} \mathcal{F}(\nabla\phi , \nabla\psi)= \int \Gamma_2(\phi, \psi) d \mu.
$$
From the work of~\cite{erbar-kuwada2015},  the curvature-dimension condition $CD(\rho,n)$ ($\rho\in\R$, $n>0$) is equivalent to the differential inequality
\begin{equation}
\label{eq-103}
\forall \mu \in \mathcal{P}_2(N), \,  \nabla\phi\in T_{\mu} \mathcal P_2(N), \quad \mathrm{Hess}_{\mu} \mathcal{F}(\nabla\phi, \nabla\phi) \geq \rho|\nabla\phi|_{\mu}^2+ \frac{1}{n} \langle \grad_{\mu} \mathcal F, \nabla\phi \rangle_{\mu}^2.
\end{equation}
%From the work of~\cite{erbar-kuwada2015}, we know that  the Manifold $(N,\mathfrak g)$ satisfies the curvature-dimension condition $CD(\rho,n)$ if and only if the entropy functional $\mathcal F$ satisfies the differential equation~\eqref{eq-103}. 
As in a  finite dimensional Riemannian manifold, the acceleration of a curve is defined as the covariant derivative of the velocity field along the curve itself. Recalling the definition of velocity $\dot{\mu}_t$ we gave at \eqref{eq-100}, it turns out that the acceleration, which we denote $\ddot{\mu}_t$ is given by 
\begin{equation}\label{eq:acc}
\ddot{\mu}_t= \nabla \left( \frac{d}{dt} \varphi_t+ \frac{1}{2} |\nabla \varphi_t|^2\right) \in T_{\mu_t}\mathcal{P}_2(N),
\end{equation}
where $\dot{\mu}_t=\nabla\varphi_t$ in the above equation.

\section{Some simple inequalities on a toy model}
\label{sec-2}

 In this section we show a series of inequalities inspired by a toy model proposed in \cite{gentil-leonard2020}. This model, in a very simple way,  captures some geometric features of the Schr\"odinger problem. Such inequalities will translate effortlessly to the true Schr\"odinger problem by means of Otto calculus, at least formally, and in this latter form they will be seen to be equivalent to local logarithmic inequalities, thus delivering the main message of this note. We shall do so in Section~\ref{sec-3}.

To build the toy model we consider a function $F: \R^d \longrightarrow \R$ be a twice differentiable function with $d>0$. We note $F'$ (resp. $F''$) the gradient (resp. the Hessian) of $F$.   For every $T>0$ and $x,y \in \R^d$, the toy model is the following optimization problem 
\begin{equation}
\label{eq-150}
C_T(x,y)=\inf \left\{ \int_0^T \SBRA{{|\dot \omega_t|^2}+{|F'(\omega_t)|^2}}dt\right\},
\end{equation}
where the infimum taken over all smooth  paths from $[0,T]$ to $\R^d$ such that $\omega_0=x$ and $\omega_T=y$. A standard variational argument shows that any minimizer  $\left( X_t^T \right)_{t \in [0,T]}$  of~\eqref{eq-150} satisfies Newton's system
\begin{equation}
\label{eq-151}
\left\{
\begin{array}{l}
\disp\ddot X_t^T =\frac{1}{2} (|F'|^2)'(X_t^T)=F''(X_t^T)F'(X_t^T),\\
\disp X_0^T=x,\,\,X_T^T=y,
\end{array}
\right.
\end{equation}
and is called an $F$-interpolation between $x$ and $y$. From \eqref{eq-151} we also deduce that 
$$
E_T(x,y)=|\dot X^T_t|^2-|F'(X_t^T)|^2,
$$
is conserved along optimal curves, i.e. it is a constant function of $t$.  $E_T$ plays the role of the conserved total (kinetic $+$ potential) energy of a physical system. Mimicking~\eqref{eq-103}, we say that $F$ is $(\rho,n)$-convex if
$$
F''\geq \rho{\rm Id}+\frac{1}{n}F' \otimes F'.
$$

We proceed to show how the announced inequalities hold under either in the $(\rho,\infty)$ or the $(0,n)$ setting.  In both cases, they will be consequences of convexity estimates for the function 
$$
\lambda(t)=\int_0^t|\dot X^T_s+F'(X_s^T)|^2ds,
$$
where  $(X_t^T)_{t\in[0,T]}$ is an $F$-interpolation. Note that $\lambda(T)=C_T(x,y)+2(F(y)-F(x))$.  Note that $\lambda(t)$ is a measure of how much the $F$-interpolation deviates from being a gradient flow.

\begin{enumerate}[\bf 1.]
\item {\bf Estimates under  $(\rho,\infty)$-convexity with $\rho\in\R$, cf.~\cite{conforti2019,gentil-leonard2020}.}
A straightforward calculation using  $(\rho,\infty)$-convexity and Newton's law yields that $\lambda$ verifies the differential inequality
$$\frac{d^2}{d t^2}\lambda\geq 2\rho\frac{d}{d t}\lambda, $$ 
which can be integrated in various ways giving the following bounds
$$
\frac{d}{d t}\lambda (0)\leq e^{-2\rho T}\frac{d}{d t}\lambda (T), \qquad  \lambda(T)-\lambda(0)\leq  \frac{1-e^{-2\rho T}}{2\rho }\frac{d}{d t} \lambda (T),$$

and
$$
\lambda(T)-\lambda(0)\geq  \frac{e^{2\rho T}-1}{2\rho }\frac{d}{d t}\lambda(0).
$$
We have therefore established the inequalities
\begin{equation}
\label{eq-152}
|\dot X_0^T+F'(x)|^2\leq e^{-2\rho T}|\dot X_T^T+F'(y)|^2,
\end{equation}
$$
C_T(x,y)+2(F(y)-F(x))\leq \frac{1-e^{-2\rho T}}{2\rho T}|\dot X_T^T+F'(y)|^2,
$$
$$
C_T(x,y)+2(F(y)-F(x))\geq \frac{e^{2\rho T}-1}{2\rho T}|\dot X_0^T+F'(x)|^2.
$$

\item {\bf Estimates under the  $(0,n)$-convexity with  $n>0$, cf.~\cite{ccg20}.}
In this case we consider the map $\Phi(t)=\lambda(t)-tE_T(x,y)$  and observe that $(0,n)$ and Newton's law combined yield
$$\Phi''\geq\frac{1}{2n}\Phi'^2,$$ 
which is equivalent to concavity of the map $[0,T]\ni t\mapsto \exp\PAR{-\frac{\Phi(t)}{2n}}$. But then, using the basic convexity estimates
$$
\exp\PAR{\frac{1}{2n}(\Phi(T)-\Phi(0))}\leq 1+\frac{1}{2n}T\Phi'(T) \quad \text{and} \qquad \exp\PAR{\frac{1}{2n}(\Phi(0)-\Phi(T))}\leq 1-\frac{1}{2n}T\Phi'(0),
$$
we obtain the inequalities
$$
\exp\PAR{\frac{1}{2n}\SBRA{C_T(x,y)+2F(y)-2F(x)-TE_T(x,y)}}\leq 1+\frac{T}{2n}\PAR{|\dot X_T^T+F'(y)|^2-E_T(x,y)},
$$
$$
\exp\PAR{-\frac{1}{2n}\SBRA{C_T(x,y)+2F(y)-2F(x)-TE_T(x,y)}}\leq 1-\frac{T}{2n}\PAR{|\dot X_0^T+F'(x)|^2-E_T(x,y)}.
$$

In particular, 
$$
E_T(x,y)\leq \frac{2n}{T}+|\dot X_0^T+F'(x)|^2,
$$
$$
E_T(x,y)\geq -\frac{2n}{T}+|\dot X_T^T+F'(y)|^2.
$$
\end{enumerate}

\section{The Schr\"odinger problem}
\label{sec-3}

We now move on to analyse the Schr\"odinger problem which we obtain from the toy model discussed above by replacing the underlying space is $\mathcal P_2(N)$ and equipping it with the Otto metric, see \cite{gigli-tamanini2020,gentil-leonard2020}. We remark that this is not the original formulation of the Schr\"odinger problem which, by means of Sanov's theorem and large deviations theory is usually cast as an entropy minimization problem on path space, see the survey \cite{leonard2014}.

\begin{edefi}[Schr\"odinger problem]
\label{def-1}
For any probability measures $\mu,\nu\in\mathcal P_2(N)$, we define te entropic cost  $\mathcal{C}_T(\mu,\nu)$ as
\begin{equation}
\label{eq-200}
 \mathcal{C}_T(\mu,\nu)=\inf\left\{\int_0^T\SBRA{{|\dot\mu_s|_{\mu_s}^2}+\mathbf{\Gamma}(\mathcal{F})(\mu_s)}ds\right\}\in[0,\infty],
\end{equation}
where the infimum runs over all absolutely continuous paths $(\mu_s)_{s\in[0,T]}$ satisfying $\mu_0=\mu$ and $\mu_T=\nu$ and   $\mathbf{\Gamma}(\mathcal{F})$ has been defined in~\eqref{eq-104}. Minimizers of $ \mathcal{C}_T(\mu,\nu)$ are called entropic interpolations between $\mu$ and $\nu$.
\end{edefi}

Let us now recall some fundamental properties of entropic interpolations. For a detailed account, we refer to~\cite{leonard2014}.
\begin{itemize}
\item For any $\mu, \nu\in \mathcal P(N)$,  absolutely continuous with a positive  $\mathcal C^\infty$ density with respect to $dx$, then there exists
 $f,g\in\mathcal C^\infty(N,(0,\infty))$ such that minimizers of~\eqref{eq-200} $(\mu_t^T)_{t\in[0,T]}$ are given by, 
\begin{equation}
\label{eq-201}
\mu_t^T=P_tfP_{T-t}g\,dx.
\end{equation}
Moreover, the couple $(f,g)$ is unique up to a multiplicative constant and the velocity field of the entropic interpolation is given by  
\begin{equation}
\label{eq-211}
\dot \mu_t^T=\nabla \log \frac{P_{T-t}g}{P_tf}.
\end{equation}

\item Entropic interpolations are solution to a second order equation, akin to Newton's law \eqref{eq-151}.  More precisely, it has been proven in~\cite[Theorem 1.2]{conforti2019} (see also~\cite[Sec. 3.3]{gentil-leonard2020} and let us notice that the acceleration of the entropic interpolation has been computed previously in~\cite{gigli-tamanini21}) that the entropic interpolation $\left( \mu_t^T \right)_{t \in [0,T]}$  is a solution of
\begin{equation} 
\label{eq-202}
\ddot \mu_t^T= \frac{1}{2} \grad_{\mu_t^T}\mathbf{\Gamma}(\mathcal{F})=\mathrm{Hess}_{\mu_t^T}\mathcal{F}\big(\grad_{\mu_t^T}\mathcal{F}\big) \, \in T_{\mu^T_t} \mathcal P_2(N),
\end{equation}
where $(\ddot \mu_t^T)_{t\in[0,T]}$ is the acceleration defined in~\eqref{eq:acc}, that is completely analogous to Newton's law~\eqref{eq-151}. Moreover, the quantity
\begin{equation} 
\label{eq-212}
\mathcal E_T(\mu,\nu)= |\dot\mu_t^T|_{\mu_t^T}^2-\mathbf{\Gamma}(\mathcal{F})(\mu_t^T)=|\dot\mu_t|_{\mu_t^T}^2-|\grad_{\mu_t^T}\mathcal F|^2_{\mu_t^T}
\end{equation}
is conserved, i.e. it doesn't depends on $t\in[0,T]$.

\end{itemize}
 It remains to write the inequalities derived for the toy model in terms of entropic interpolations. By replacing the function $\lambda(t)$ used in section 2 with the function 
$$
\Lambda(t)=\int_0^t|\dot \mu_s^T+\grad_{\mu_s^T} \mathcal F|^2_{\mu_s^T}ds,
$$
where  $({\mu_s^T})_{t\in[0,T]}$ is an entropic interpolation. We obtain, arguing exactly as before but working on $\mathcal{P}_2(N)$ and using the Riemannian formalism associated with the Otto metric the following estimates %It satisfies  $\Lambda(T)=\mathcal C_T(x,y)+2(\mathcal F(\nu)-\mathcal F(\mu))$.

\medskip

%Under the under the  $CD(\rho,\infty)$ condition, $\Lambda$ satisfies $\Lambda''\geq 2\rho\Lambda'$ (which has been noticed in~\cite{conforti2019}) and under the  $CD(0,n)$ condition the maps $\Phi(t)=\Lambda(t)-t\mathcal E_T(\mu,\nu)$ satisfies $\Phi''\geq\frac{1}{2n}\Phi'^2$ (which has been noticed in~\cite{ccg20}). Using again Lemma~\ref{lem-1}, we get the two next results.  
\begin{enumerate}[1.]
\item {\bf Estimates under  $CD(\rho,\infty)$ with $\rho\in\R$, cf.~\cite{conforti2019,gentil-leonard2020}.} In this case we have that if $(\mu_t^T)_{t\in[0,T]}$ is the entropic interpolation between $\mu$ and $\nu$:

\begin{equation}
\label{eq-203}
|\dot \mu_0^T+\grad_{\mu}\mathcal F|^2_\mu\leq e^{-2\rho T}|\dot \mu_T^T+\grad_{\nu}\mathcal F|^2_\nu,
\end{equation}
\begin{equation}
\label{eq-204}
\mathcal C_T(\mu,\nu)+2(\mathcal F(\nu)-\mathcal F(\mu))\leq \frac{1-e^{-2\rho T}}{2\rho }|\dot \mu_T^T+\grad_{\nu}\mathcal F|^2_\nu,
\end{equation}
\begin{equation}
\label{eq-205}
\mathcal C_T(\mu,\nu)+2(\mathcal F(\nu)-\mathcal F(\mu))\geq \frac{e^{2\rho T}-1}{2\rho }|\dot \mu_0^T+\grad_{\mu}\mathcal F|^2_\mu.
\end{equation}

\item {\bf Estimates under the  $(0,n)$-convexity with  $n>0$, cf.~\cite{ccg20}.}  In this case we have that if $(\mu_t^T)_{t\in[0,T]}$ is the entropic interpolation between $\mu$ and $\nu$:
\begin{equation}
\label{eq-206}
\exp\PAR{\frac{1}{2n}\SBRA{\mathcal C_T(\mu,\nu)+2\mathcal F(\nu)-2\mathcal F(\mu)-T\mathcal E_T(\mu,\nu)}}\leq 1+\frac{T}{2n}\PAR{|\dot \mu_T^T+\grad_{\nu}\mathcal F|^2_\nu-\mathcal E_T(\mu,\nu)},
\end{equation}
\begin{equation}
\label{eq-207}
\exp\PAR{-\frac{1}{2n}\SBRA{\mathcal C_T(\mu,\nu)+2\mathcal F(\nu)-2\mathcal F(\mu)-T\mathcal E_T(\mu,\nu)}}\leq 1-\frac{T}{2n}\PAR{|\dot \mu_0^T+\grad_{\mu}\mathcal F|^2_\mu-\mathcal E_T(\mu,\nu)}.
\end{equation}

In particular, 
$$
\mathcal E_T(\mu,\nu)\leq \frac{2n}{T}+|\dot \mu_T^T+\grad_{\nu}\mathcal F|^2_\nu,
$$
\begin{equation}\label{eq-variationalLiYau}
\mathcal E_T(\mu,\nu)\geq -\frac{2n}{T}+|\dot \mu_0^T+\grad_{\mu}\mathcal F|^2_\mu.
\end{equation}
\end{enumerate}

\section{Local Logarithmic Sobolev inequalities via Schr\"odinger problem and Otto calculus}
\label{sec-4}
The main contribution of this note is to show how the Bakry-\'Emery-Ledoux estimates~\eqref{eq-10} to~\eqref{eq-70} are in fact equivalent to the inequalities derives in Section~\ref{sec-3} along entropic interpolations. The key to do this is the product formula \eqref{eq-201}.

%For any probability measures $\mu,\nu$, absolutely continuous with a smooth densities with respect to the measure $dx$, the entropic interpolation is given by 
%$$
%\mu_t^T=P_t f P_{T-t}g\,dx,\,\,t\in[0,T],
%$$ 
%for some functions $f,g>0$. We recall  that $\dot \mu_t^T=\nabla \log\frac{P_{T-t}g}{P_t f }$ and $\grad_{\mu_t^T}\mathcal F=\nabla\log({P_{T-t}g}{P_t f })$. 

%\medskip

\paragraph{The mild gradient commutation estimate \eqref{eq-10}.} 
Let us investigate in details the first inequality~\eqref{eq-203}.  Using \eqref{eq-201} we find  
\begin{equation}
\label{eq-312}
| \dot\mu_t^T+\grad_{\mu_t^T}\mathcal F|^2_{\mu_t^T}\!=\!\int\! \Gamma\PAR{\log\frac{P_{T-t}g}{P_tf}+\log({P_{T-t}g}{P_tf})}{P_{T-t}g}{P_tf}\,dx\!=\!4\!\int\! \frac{\Gamma\PAR{{P_{T-t}g}}}{P_{T-t}g}{P_tf}\,dx.
\end{equation}
Then, by reversibility of the semigroup $(P_t)_{t\geq 0}$,  inequality~\eqref{eq-203} can be written as 
\begin{equation}
\label{eq-322}
\int\frac{\Gamma( P_Tg)}{P_Tg} f\,dx\leq e^{-2\rho T}\int P_T\Big(\frac{\Gamma( g)}{g}\Big)f \,dx.
\end{equation}
But, for any positive functions $f,g$, the path 
\begin{equation}
\label{eq-332}
[0,T]\ni t\mapsto \frac{P_tfP_{T-t}g\,dx}{\int\! fP_Tg \,dx}
\end{equation}
is an interpolation entropic, then the previous inequality~\eqref{eq-322} is valid for any functions $f,g>0$.
If we now let  $f\,dy \rightarrow \delta_x$ weakly, we get
$$
\frac{\Gamma( P_Tg)}{P_Tg}(y)\leq e^{-2\rho T}P_T\Big(\frac{\Gamma( g)}{g}\Big)(y),
$$ 
which is precisely \eqref{eq-10}.

\paragraph{The local logarithmic Sobolev inequalities \eqref{eq-20}-\eqref{eq-30} under $CD(\rho,\infty)$. } 

Let us have a look to the inequality~\eqref{eq-204}. First, since $(\mu_t^T)_{t\in[0,T]}$ is the entropic interpolation between $\mu$ and $\nu$, then 
$$
\mathcal C_T(\mu,\nu)+2(\mathcal F(\nu)-\mathcal F(\mu))=\int_0^T| \dot\mu_s^T+\grad_{\mu_s^T}\mathcal F|^2_{\mu_s^T}ds,
$$
and from~\eqref{eq-312}, 
$$
\int_0^T| \dot\mu_t^T+\grad_{\mu_t^T}\mathcal F|^2_{\mu_t^T}dt=
4\int_0^T\int \frac{\Gamma\PAR{P_{T-t}g}}{P_{T-t}g}{P_tf}\,dxdt=4\int_0^T\int P_t\PAR{\frac{\Gamma\PAR{P_{T-t}g}}{P_{T-t}g}}f\,dxdt.
$$
We know that if $\Phi(t)=P_t(P_{T-t}g\log P_{T-t}g)$, ($t\in[0,T]$), then  
$$
\Phi'(t)=P_t\PAR{\frac{\Gamma\PAR{P_{T-t}g}}{P_{T-t}g}}.
$$
So, at the end,  
\begin{equation}
\label{eq-302}
\mathcal C_T(\mu,\nu)+2(\mathcal F(\nu)-\mathcal F(\mu))=4\int\int_0^T \Phi'(t)dtf\,dx=4\int [P_T(g\log g)-P_Tg\log P_Tg]f\,dx.
\end{equation}
By using the path defined in~\eqref{eq-332}, one can conclude that inequality~\eqref{eq-204} can be written as
$$
\int \SBRA{P_T(g\log g)-P_Tg\log P_Tg}f\,dx\leq \frac{1-e^{-2\rho T}}{2\rho }\int P_T\PAR{\frac{\Gamma(g)}{g}}f\,dx, 
$$
for any functions $f,g>0$. Again, if $f\,dx\rightarrow\delta_y$ ($y\in N$) weakly, the previous inequality becomes 
$$
P_T(g\log g)(y)-P_Tg(y)\log P_Tg(y)\leq \frac{1-e^{-2\rho T}}{2\rho } P_T\PAR{\frac{\Gamma(g)}{g}}(y),
$$
which is exactly the local logarithmic Sobolev inequality~\eqref{eq-20}. Arguing similarly, inequality~\eqref{eq-205} is seen to be equivalent to the local reverse logarithmic Sobolev inequality~\eqref{eq-30}.   

\medskip

\noindent
\paragraph{Local logarithmic Sobolev inequalities~\eqref{eq-40}-\eqref{eq-60} under $CD(0,n)$ and the Li-Yau inequality~\eqref{eq-70}.} 

To understand the last inequalities, is it enough to compute the conserved quantity. We have, 
\begin{multline}
\label{eq-303}
\mathcal E_T(\mu,\nu)=\int\SBRA{\Gamma\PAR{\log\frac{P_{T-t}g}{P_tf}}-\Gamma\PAR{\log{P_{T-t}g}{P_tf}}}{P_{T-t}g}{P_tf}\,dx=\\-4\int\Gamma(\log{P_{T-t}g},\log{P_tf}){P_{T-t}g}{P_tf}\,dx=-4\int\Gamma(P_{T-t}g,P_tf)\,dx=4\int\Delta P_{T}gf\,dx.
\end{multline}
From the identities~\eqref{eq-302}, \eqref{eq-312} and~\eqref{eq-303}, inequalities~\eqref{eq-206} becomes, 
\begin{multline*}
\exp\PAR{\frac{1}{2n}\SBRA{4\int  \SBRA{P_T(g\log g)-P_Tg\log P_Tg}f\,dx-4T\int\Delta P_{T}gf\,dx}}\\
\leq 1+\frac{T}{2n}\PAR{4\int P_T\Big(\frac{\Gamma( g)}{g}\Big)f \,dx-4\int\Delta P_{T}gf\,dx}.
\end{multline*}
Using the path of probability measures~\eqref{eq-332}, for any functions $f,g>0$, we have
\begin{multline*}
\exp\PAR{\frac{2}{n\int P_Tg f\,dx}\SBRA{\int  \SBRA{P_T(g\log g)-P_Tg\log P_Tg}f\,dx-T\int\Delta P_{T}gf\,dx}}\\
\leq 1+\frac{2T}{n\int P_Tg f\,dx}\PAR{\int P_T\Big(\frac{\Gamma( g)}{g}\Big)f \,dx-\int\Delta P_{T}gf\,dx}.
\end{multline*}
Again, if we let $fdx\rightarrow\delta_y$ ($y\in N$), we get inequality~\eqref{eq-40}. A similar argument can be used for the equivalence between \eqref{eq-207} and the reverse local inequality~\eqref{eq-50} and to show the equivalence of \eqref{eq-variationalLiYau} and the Li-Yau inequality \eqref{eq-70}.

\begin{erem}
\label{rem-2}
When $\mu=\delta_y$ for some $y\in N$ and $\nu\in P(N)$, the entropic interpolation is explicit and given by the following path 
\begin{equation}
\label{eq-300}
[0,T]\ni t\mapsto \mu_t^T=p_t^y P_{T-t}\Big(\frac{d\nu}{dx}\frac{1}{p_T^y}\Big)\,dx.
\end{equation}
It is interesting to notice that it is in fact, the path used in the computations of this note.

%When we choose $f\,dx=\delta_y$ in the previous computations, the entropic interpolation becomes
%$$
%[0,T]\ni t\mapsto \mu_t^T=p_t^y P_{T-t}g\,dx.
%$$
%This choice implies that  $\mu_0^T=\delta_y=\mu$ and then we have $g=\frac{d\nu}{\,dx}\frac{1}{P_T^y}$.  The path
%\begin{equation}
%\label{eq-300}
%[0,T]\ni t\mapsto \mu_t^T=p_t^y P_{T-t}\Big(\frac{d\nu}{dy}\frac{1}{p_T^y}\Big)\,dx,
%\end{equation}
%is the entropic interpolation between between $\delta_y$ and $\nu$. Functions $f$ and $g$ are explicit! In the case where the final measure is also a Dirac measure at some $z\in N$, then  
%\begin{equation}
%\label{eq-301}
%[0,T]\ni t\mapsto \mu_t^T=\frac{p_t^y p_{T-t}^z}{p_T^y(z)}\,dx,
%\end{equation}
%is the entropic interpolation between $\delta_y$ and $\delta_z$. 
%Of course, for these interpolations the entropic cost function $\mathcal C_T$ are infinite since the entropy is infinite at a Dirac measure. One way to deals with a finite cost function is to modify it by defining  the relaxed cost function 
%$$
%\inf\BRA{\int |\dot\mu_s^T+\grad_{\mu_s^T}\mathcal F|^2_{\mu_s^T}ds}, 
%$$
%where the infimum runs over all absolutely continuous paths $(\mu_s^T)_{s\in[0,T]}$ from $\mu$ to $\nu$. This minimizer problem is the same as the usual one, up to entropic term as it is explained for instance in~\cite{gentil-leonard2020}. In particular, the entropy of the initial or the final values don't need to be finite. 
\end{erem}

\paragraph{Conclusion}
We have demonstrated how local logartihmic Sobolev inequalities~\eqref{eq-10}-\eqref{eq-50} admit a nice geometrical interpretation via inequalities~\eqref{eq-203}-\eqref{eq-207} by means of Otto calculus and the Schr\"odinger problem. Such interpretation was first given for the integrated version of the logarithmic Sobolev inequality in \cite{jko1998} and our results are the first, to the best of our knowledge to cover the local versions. We hope that our considerations will be of help to understand local inequalities for general functionals satisfying the curvature dimension conditions beyond, thus going beyond the case of the Boltzmann entropy.

\paragraph{Acknowledgements} This reflexion  was supported by the French ANR-17-CE40-0030 EFI project.

\end{document}